\documentclass[12pt,a4paper]{article}
\usepackage[utf8]{inputenc}
\usepackage{tikz}
\usepackage{amsfonts}
\usepackage{float,epsfig,floatflt,here}
\usepackage{amsmath}
\usepackage{amssymb}
\usepackage{layout}
\usepackage{graphicx}
\usepackage{hyperref}
\usepackage[a4paper, total={6in, 8in}]{geometry}
\usepackage{mathtools}
\usepackage[T1]{fontenc}
\usepackage{babel}
\usepackage{multirow}
\usepackage{array}
\newtheorem{definition}{Definition}[section]
\newtheorem{theorem}{Theorem}
\newtheorem{example}{Example}

\DeclareMathOperator{\spn}{span}
\title{\textbf{Dynamical Representation of Frames in Tensor Product of Hardy Spaces }}

\begin{document}

\maketitle

\author{  
\begin{tabular}{*{2}{>{\centering}p{.5\textwidth}}}
\large Nabin Kumar Sahu$^*$  & \large Vishesh Kumar$^*$   \tabularnewline
\href{mailto: nabinkumar\_sahu@daiict.ac.in}{nabinkumar\_sahu@daiict.ac.in} & \href{mailto:visheshrajput771@gmail.com}{visheshrajput771@gmail.com} 
\end{tabular}}\\

$^*${\centering{Dhirubhai Ambani Institute of Information and Communication Technology, Gandhinagar, India}}\\

\def\thefootnote{*}\footnotetext{Vishesh Kumar and Nabin Kumar Sahu are co-first authors with equal contribution.}\def\thefootnote{\arabic{footnote}}

\begin{abstract}
Dynamical Sampling of frames and tensor products are important topics in harmonic analysis. This paper combines the concepts of dynamical sampling of frames and the Carleson condition in the tensor product of Hardy spaces. Initially we discuss the preservation of the frame property under the tensor product on the Hilbert spaces. Then we discuss the iterative representation of frames in tensor product of Hardy spaces. The key ingredient of this paper is the so-called Carleson condition on the sequence $\{ \lambda_k \}_{k=1}^{\infty} \otimes\{ \gamma_l \}_{l=1}^{\infty} $ in the open unit disc $\mathbb{D}_1 \otimes \mathbb{D}_2$. Our proof is motivated by the result of Shapiro and Shields.
\end{abstract}



\section{Introduction}

Dynamical sampling of frames is a technique used in signal processing that deals with recovering a signal from samples taken in space and time. This technique aims to accurately represent the evolution process of the signal over time by selecting frames that best capture its dynamic behavior. By capturing the dynamic behavior of the signal over time and space, dynamical sampling of frames can help improve the quality of signal recovery in scenarios where an evolution process generates the signal. It is a relatively new topic in applied harmonic analysis but has already garnered considerable interest among researchers and practitioners alike \cite{aldroubi2017dynamical, aldroubi2017dynamical1, aldroubi2013dynamical,cabrelli2020dynamical, carando2009duality,martin2021continuous,nemec2017quantitative,stoeva2011characterization}. \\

\noindent Aldroubi et al. \cite{aldroubi2017dynamical}  initiated the study of dynamical sampling problem and the frame properties of sequences in the form of $\{ T^n \varphi \}_{n=0}^{\infty},$ where $T: \mathcal{H} \to \mathcal{H} $ belongs to certain classes of linear operators and $\varphi \in \mathcal{H}$, where $\mathcal{H}$ is a separable Hilbert space. The classical dynamical sampling problem is as follows: Consider a measurement set $Y=\{x(i),Tx(i),T^{2}x(i),....., T^{m_{i}}x(i):~i\in \Omega\}$, where $T$ is a bounded operator on $l^{2}(I)$, $\Omega$ is an index set. The classical dynamical sampling problem is to find the necessary and sufficient conditions on $A$, $\Omega$ and $m_{i}$ such that any $f\in l^{2}(I)$ can be recovered efficiently from the measurement set $Y$.

In 2017,  Christensen and Hasannasab \cite{christensen2017operator} classified the frames which is generated by a linear operator not necessarily bounded. They have also discussed dynamically generated dual frames and they found that the dynamical frames are unstable under the classical perturbation conditions (Paley-Wiener type conditions).
Again in 2019, Christensen and Hasannasab \cite{christensen2019frame} proved that every frame
which is norm-bounded below can be represented as a finite union of sequences $\{(T_{j})^{n}\phi_{j}\}_{n=0}^{\infty}$ for some bounded operators $T_{j}$ and some elements $\phi_{j}$
  in the underlying Hilbert space. 
Recently,  Ashbrock and  Powell \cite{ashbrock2023dynamical} proved that every redundant finite frame for $\mathbb{F}^d$, where $\mathbb{F} = \mathbb{R}$ or $\mathbb{C},$ has infinitely many dynamical dual frames, and introduced a low complexity error diffusion quantization algorithm based on dynamical dual frames.\\

Hardy Hilbert space $H^{2}(\mathbb{D})$ is the space of all analytic functions defined on the unit disk $\mathbb{D}$ with square summable coefficients. That is $H^{2}(\mathbb{D})= \Big\{f:\mathbb{D}\rightarrow \mathbb{C}:~~f(z)=\sum_{n=0}^{\infty}a_{n}z^{n},~\text{and}~\sum_{n=0}^{\infty}|a_{n}^{2} |< \infty\Big\}$. The space $H^{2}(\mathbb{D})$ is equipped with the inner product $\displaystyle\langle f(z),g(z)\rangle=\sum_{n=0}^{\infty}a_{n}\overline{b_{n}}$, where $f(z)=\displaystyle \sum_{n=0}^{\infty}a_{n}z^{n}$ and $g(z)=\displaystyle\sum_{n=0}^{\infty}b_{n}z^{n}$ .
The theory of frames has been well investigated in the space $H^{2}(\mathbb{D})$, and for that one may refer \cite{chen2017frame, speransky2019existence, volkmer1995frames} and the references there in. Carleson condition \cite{christensen2018operator} is the most important ingredient to construct dynamically generated frames in a separable Hilbert space. We recall that a sequence $\displaystyle\{\alpha_{k}\}_{k=1}^{\infty}\subset \mathbb{D}$ satisfies the Carleson condition if $$\displaystyle \inf_{n\in \mathbb{N}}\displaystyle \Pi_{k\neq n}\frac{|\alpha_{k}-\alpha_{n}|}{|1-\overline{\alpha_{k}}\alpha_{n}|}>0.$$ We quote the main result from  \cite{christensen2018operator}.
\begin{theorem} \rm Let $\mathcal{H}$ be a separable Hilbert space with orthonormal basis $\displaystyle \{e_{k}\}_{k=1}^{\infty}$. 
Let $\displaystyle\{\alpha_{k}\}_{k=1}^{\infty}$ be a sequence in the unit disk $\mathbb{D}$ in the complex plane. Assume that $\displaystyle \big\{\sqrt{1-|\alpha_{k}|^{2}}\big\}_{k=1}^{\infty}\in l^{2}(\mathbb{N})$. Let $T:\mathcal{H}\rightarrow \mathcal{H}$ be a bounded linear operator such that $Te_{k}=\alpha_{k}e_{k}$. Then the sequence $\displaystyle \{T^{n}h\}_{n=0}^{\infty}=\Big\{ \sum_{k=1}^{\infty} \alpha_{k}^{n}\sqrt{1-|\alpha_{k}|^{2}}e_{k}\Big\}$ is a frame for the Hilbert space $\mathcal{H}$ if and only if $\displaystyle\{\alpha_{k}\}_{k=1}^{\infty}$ satisfies the Carleson condition.
\end{theorem}

\noindent In \cite{bourouihiya2008tensor}, it is proved that the tensor product of two sequences is frame if and only if each part of this product is a frame. In \cite{khosravi2012frames}, it is also proved that if $\{x_n \}_{n \in I}$ and $\{y_m \}_{m \in I} $ be a frame for $\mathcal{H}_1$ and $\mathcal{H}_2$, respectively. Then $\{ x_n \otimes y_m \}_{n,m \in I}$ is a frame for $\mathcal{H}_1 \otimes \mathcal{H}_2.$

\par
The above work persuade us to carry out  investigation on dynamically generated frames on tensor product of two Hilbert spaces. The main contribution of this report is Section 3, where the Carleson condition has been newly formulated. The relationship between the Carleson condition and the frame properties of an iterated sequence has been established.


\section{Iterative representation of sequences in tensor product of two Hilbert spaces}
In this section we consider linearly  independent sequences $\{f_k \}_{k=1}^{\infty} \otimes \{ g_k \}_{k=1}^{\infty}$ in Hilbert space $\mathcal{H}_1 \otimes \mathcal{H}_2$ and explore the existence of a representation $\{ T_{1}^{n} f_1 \}_{n=0}^{\infty} \otimes \{T_{2}^{n} g_1 \}_{n=0}^{\infty} $ with a bounded operator $T_1 \otimes T_2 : span \{f_k \}_{k=1}^{\infty} \otimes span \{g_k \}_{k=1}^{\infty} \to span \{f_k \}_{k=1}^{\infty} \otimes span \{g_k \}_{k=1}^{\infty}.$ This generalization is motivated by the results for frames proved in \cite{christensen2018operator,christensen2018dynamical}. In general the existence of a representation of a sequence $\{ f_k \}_{k=1}^{\infty} \otimes \{ g_k \}_{k=1}^{\infty} $ of the form $\{T_{1}^n f_1 \}_{n=0}^{\infty} \otimes \{ T_{2}^n g_1 \}_{n=0}^{\infty} $ for a bounded operator $ T_1 \otimes T_2$ is not nearly related with frame properties of the given sequence $ \{f_k \}_{k=1}^{\infty} \otimes \{ g_k \}_{k=1}^{\infty}.$

\begin{example}\rm
Consider a linearly independent frame $ \{f_k \}_{k=1}^{\infty} \otimes \{ g_k \}_{k=1}^{\infty}$ for an infinite dimensional Hilbert space $\mathcal{H}_1 \otimes \mathcal{H}_2$ and the associated representation $ \{ f_k \}_{k=1}^{\infty} \otimes \{g_k \}_{k=1}^{\infty} = \{ T_{1}^{n} f_1 \}_{n=0}^{\infty} \otimes \{ T_{2}^n g_1 \}_{n=0}^{\infty} $ in terms of a linear operator $ T_1 \otimes T_2 : \spn \{f_k \}_{k=1}^{\infty} \otimes \{g_k \}_{k=1}^{\infty} \rightarrow \spn \{f_k \}_{k=1}^{\infty} \otimes \{g_k \}_{k=1}^{\infty} $. Assume that $\inf_{k\in \mathbb{N}} || f_k \otimes g_k || =  \inf_{k \in \mathbb{N}} ||f_k||||g_k|| > 0$. Consider the sequence $\{ \phi_k \}_{k=1}^{\infty} \otimes \{ \psi_{k} \}_{k=1}^{\infty} \subset \mathcal{H}_1 \otimes \mathcal{H}_2$ given by $\phi_k \otimes \psi_k = 2^k f_k \otimes 2^k g_k $ $, k \in \mathbb{N},$ which leads to a frame like expansion and satisfies the lower frame condition but fails the upper one. For any $k \in \mathbb{N},$
 \begin{align*}
    \phi_{k+1} \otimes \psi_{k+1} &= 2^{k+1} f_{k+1} \otimes 2^{k+1} g_{k+1}\\
    &= 2^{k+1}T_1 f_k \otimes 2^{k+1} T_2 g_k\\
    &=4 T_1 (2^k f_k) \otimes T_2 (2^k g_k)\\
    &= 4 T_1 ( \phi_k ) \otimes T_2 (\psi_k )\\
    &= 4 (T_1 \otimes T_2 )( \phi_k \otimes \psi_k ).\\
\end{align*}

This shows that $\{ \phi_k \}_{k=1}^{\infty} \otimes \{ \psi_{k} \}_{k=1}^{\infty} $ has the representation  $\{ W_{1}^n \phi_1 \}_{n=0}^{\infty} \otimes \{ W_{2}^n \psi_1 \}_{n=0}^{\infty}$, where $W_1 = 2T_1, W_2 = 2T_2.$ In particular we can say that the frame $\{ f_k \}_{k=1}^{\infty} \otimes \{ g_k \}_{k=1}^{\infty}$ is represented by a bounded operator if and only if the non-Bessel sequence $ \{ \phi_k \}_{k=1}^{\infty} \otimes \{ \psi_k \}_{k=1}^{\infty} $ is represented by a bounded operator. 
\end{example}

\par

 It is important to note a sequence $\{f_k \}_{k=1}^{\infty} \otimes \{g_k \}_{k=1}^{\infty}$ that has nice frame properties and a nice representation $\{f_k \}_{k=1}^{\infty} \otimes \{g_k \}_{k=1}^{\infty} = \{ T_{1}^n f_1 \}_{n=0}^{\infty} \otimes \{ T_{2}^{n} g_1 \}_{n=0}^{\infty}$ do not necessarily correlate with one another, as illustrated by the following example. Although the family of functions $\{f_k\}_{k=1}^{\infty} \otimes \{g_k \}_{k=1}^{\infty}$ is not a frame and cannot provide a frame-like expansion, but it has a representation $\{f_k \}_{k=1}^{\infty} \otimes \{g_k \}_{k=1}^{\infty} = \{T_{1}^n f_1 \}_{n=0}^{\infty} \otimes \{T_{2}^{n} g_1 \}_{n=0}^{\infty}$ for a bounded isometric operator $T_1 \otimes T_2.$

\begin{example}\rm
    Consider the sequence $\{f_k\}_{k=1}^{\infty} \otimes \{g_k \}_{k=1}^{\infty} $ given by $f_k \otimes g_k = (e_k + e_{k+1} ) \otimes (d_k +d_{k+1}), k \in \mathbb{N}$, where $\{e_k \}_{k=1}^{\infty} \otimes \{ d_k \}_{k=1}^{\infty}$ is an orthonormal basis for the Hilbert space $\mathcal{H}_1 \otimes \mathcal{H}_2$. We find that $\{ f_k \}_{k=1}^{\infty} \otimes \{g_k \}_{k=1}^{\infty}$ is a Bessel sequence but not a frame, despite the fact that $\overline{span}\{f_k \}_{k=1}^{\infty} \otimes \overline{span}\{g_k \}_{k=1}^{\infty} = \mathcal{H}_1 \otimes \mathcal{H}_2.$ Since $\{ f_k \}_{k=1}^{\infty} \otimes \{g_k \}_{k=1}^{\infty}$ is linearly independent, we define the operator $T_1 \otimes T_2 : span \{f_k \}_{k=1}^{\infty} \otimes span\{g_k \}_{k=1}^{\infty} \to span \{f_k \}_{k=1}^{\infty} \otimes span\{g_k \}_{k=1}^{\infty} $, by  $T_1 f_k \otimes T_2 g_k = f_{k+1} \otimes g_{k+1} ,$ and we have $\{f_k \}_{k=1}^{\infty} \otimes \{g_k \}_{k=1}^{\infty} = \{ T_{1}^n f_1 \}_{n=0}^{\infty} \otimes  \{ T_{2}^n g_1 \}_{n=0}^{\infty}.$ Then, for any $a_k , c_k \in \mathbb{C},$ and any $N \in \mathbb{N},$\\
    \begin{align*}
        \footnotesize\left|\left| T_1 \sum_{k=1}^{N} c_k f_k \otimes T_2 \sum_{k=1}^{N} a_k g_k  \right| \right|^2 &= \left| \left| (T_1 \otimes T_2)\left( \sum_{k=1}^{N} c_k f_k \otimes \sum_{k=1}^{N} a_k g_k \right) \right| \right|^2\\ 
        &= \left| \left| \sum_{k=1}^{N} c_k( e_{k+1} + e_{k+2}) \otimes a_k( d_{k+1} + d_{k+2})\right| \right|^2\\
        \displaybreak &= \left| \left| \sum_{k=1}^{N} c_k(e_k + e_{k+1}) \otimes a_k (d_k +d_{k+1})\right| \right|^2\\
        &= \left| \left| \sum_{k=1}^{N} c_k f_k \otimes \sum_{k=1}^{N} a_k g_k \right| \right|^2\\ 
        &= \left| \left|\sum_{k=1}^{N} c_k f_k\right| \right|^2 \left| \left| \sum_{k=1}^{N} a_k g_k \right| \right|^2
    \end{align*}
   It follows that $T_1 \otimes T_2$ has an extension to an isometric operator $T_1 \otimes T_2 : \mathcal{H}_1 \otimes \mathcal{H}_2 \to \mathcal{H}_1 \otimes \mathcal{H}_2.$

\end{example}

Given any sequence $\{ f_k \}_{k=1}^{\infty} \otimes \{g_k \}_{k=1}^{\infty} \subset \mathcal{H}_1 \otimes \mathcal{H}_2,$ then the synthesis operator is defined as 
\[ V : \mathcal{D}(V) \to \mathcal{H}_1 \otimes \mathcal{H}_2, V( \{ c_k \}_{k=1}^{\infty} \otimes \{d_k \}_{k=1}^{\infty}) := \sum_{k=1}^{\infty} c_k f_k \otimes \sum_{k=1}^{\infty} d_k g_k,\]
where the domain $\mathcal{D}(V) $ consists the set of all scalar-valued sequences $\{c_k \}_{k=1}^{\infty} \otimes \{d_k \}_{k=1}^{\infty}$ for which $\displaystyle \sum_{k=1}^{\infty} c_k f_k \otimes \sum_{k=1}^{\infty} d_k g_k$ is convergent.Here we do not restrict our attention to sequences $\{c_k \}_{k=1}^{\infty} \otimes \{d_k \}_{k=1}^{\infty}$ belonging to $l^2(\mathbb{N}) \otimes l^2 ( \mathbb{N} )$.\\ 

\noindent Consider the right-shift operators $\mathcal{T}_1$ and $\mathcal{T}_2$, which acts on arbitrary scalar sequences $\{c_k \}_{k=1}^{\infty}$ and $ \{d_k \}_{k=1}^{\infty}$ respectively by   $\mathcal{T}_1\{c_k \}_{k=1}^{\infty} = \{0,c_1,c_2,c_3,.......\}$ and $\mathcal{T}_2\{d_k \}_{k=1}^{\infty} = \{ 0, d_1,d_2,d_3,......\}$.Now define an operator $\mathcal{T}_1 \{c_k\}_{k=1}^{\infty} \otimes \mathcal{T}_2 \{d_k \}_{k=1}^{\infty} = \{ 0, c_1 , c_2,....\} \otimes \{ 0, d_1, d_2,....\} $. A  vector space $W$ of scalar-valued sequences $\{c_k \}_{k=1}^{\infty} \otimes \{d_k \}_{k=1}^{\infty}$ is said to be invariant under the defined operator $\mathcal{T}_1 \otimes \mathcal{T}_2$ if $(\mathcal{T}_1 \otimes \mathcal{T}_2)(W) \subseteq W$.\\

The following result generalizes one of the important results in \cite{christensen2018dynamical,christensen2018operator} for any arbitrary sequence (not necessarily frame).

\begin{theorem}\rm
    Consider a sequence $\{ f_k \}_{k=1}^{\infty} \otimes \{g_k \}_{k=1}^{\infty} \in \mathcal{H}_1 \otimes \mathcal{H}_2$ which has representation $\{ T_{1}^n f_1 \}_{n=0}^{\infty} \otimes \{ T_{2}^{n} g_1 \}_{n=0}^{\infty}$ for a linear operator $T_1 \otimes T_2 : span \{f_k\}_{k=1}^{\infty} \otimes span \{g_k \}_{k=1}^{\infty} \to span \{f_k\}_{k=1}^{\infty} \otimes span \{g_k \}_{k=1}^{\infty}.$ Then the following results hold.
\begin{enumerate}
    \item If $T_1 \otimes T_2$ is bounded, then the domain $\mathcal{D}(V)$ and the kernel $N(V)$ of the synthesis operator are invariant under the defined operator $\mathcal{T}_1 \otimes \mathcal{T}_2$, and $\left\{ \frac{|| f_{k+1} \otimes g_{k+1}||}{|| f_k \otimes g_k||} \right\}_{k=1}^{\infty} \in l^{\infty}$ where $f_k \otimes g_k \neq 0, \forall k \in \mathbb{N}.$
    \item $T_1 \otimes T_2$ is bounded on $ span \{f_k \}_{k=1}^{\infty} \otimes span \{g_k \}_{k=1}^{\infty}$ if and only if there exists a positive constant $K$ such that 
    \\ $\displaystyle \left| \left| V( \mathcal{T}_1 \{c_k \}_{k=1}^{\infty} \otimes \mathcal{T}_2 \{d_k \}_{k=1}^{\infty} ) \right| \right| \leq K \left| \left| V ( \{c_k \}_{k=1}^{\infty} \otimes 
    \{ d_k \}_{k=1}^{\infty} ) \right| \right|$ for all finite sequences $\{ c_k \}_{k=1}^{\infty} \otimes \{d_k \}_{k=1}^{\infty}.$
\end{enumerate}

\end{theorem}
\textbf{Proof:}
    We know if $T_1$ and  $T_2$ is bounded then $T_1 \otimes T_2$ is bounded, let $\widetilde{T_{1}} \otimes \widetilde{T_2}$ is a unique extension of the  bounded linear operator $T_1 \otimes T_2$ on $\overline{span}\{f_k \}_{k=1}^{\infty} \otimes \overline{span}\{g_k \}_{k=1}^{\infty}.$ \\
    
    \noindent (1) Assume $T_1 \otimes T_2$ is bounded and consider sequence $\{c_k \}_{k=1}^{\infty} \otimes \{d_k \}_{k=1}^{\infty} \in \mathcal{D}(V).$ In order to show that $\mathcal{T}_1 \{c_k \}_{k=1}^{\infty} \otimes \mathcal{T}_2 \{ d_k \}_{k=1}^{\infty} \in \mathcal{D}(V),$ i.e., that $\displaystyle \sum_{k=1}^{\infty}c_k f_{k+1} \otimes \sum_{k=1}^{\infty}d_k g_{k+1}$ is convergent. Consider any $M,N \in \mathbb{N}$ where $N > M.$ Then
\footnotesize{\begin{align*}
        \displaystyle  \left| \left| \sum_{k=1}^{N} c_k f_{k+1} \otimes \sum_{k=1}^{N}d_{k} g_{k+1} - \sum_{k=1}^{M}c_k f_{k+1} \otimes \sum_{k=1}^{M} d_k g_{k+1} \right| \right| &= \left| \left| \sum_{k=M+1}^{N} c_k f_{k+1} \otimes \sum_{k=M+1}^{N} d_k g_{k+1} \right| \right|\\
        \displaystyle &= \left| \left| T_1 \sum_{k=M+1}^{N} c_k f_k \otimes T_2 \sum_{k=M+1}^{N} d_k g_k \right| \right|\\
       \displaystyle  &=\left| \left| (T_1 \otimes T_2 ) \left( \sum_{k=M+1}^{N} c_k f_k \otimes \sum_{k=M+1}^{N} d_k g_k \right) \right| \right|\\
        &\leq \left| \left| T_1 \otimes T_2 \right| \right| \left| \left| \sum_{k=M+1}^{N} c_k f_k \right| \right| \left| \left| \sum_{k=M+1}^{N} d_k g_k \right| \right|\\
        & \longrightarrow 0
    \end{align*} } 
as $M,N \to \infty.$ Thus $ \displaystyle \sum_{k=1}^{\infty} c_k f_{k+1} \otimes \sum_{k=1}^{\infty} d_k g_{k+1}$ is convergent, i.e, $\mathcal{D}(V)$ is indeed invariant under the operator $ \mathcal{T}_1 \otimes \mathcal{T}_2$.\par
To prove the invariance of $\mathcal{N}(V),$ assume that $\{ c_k \}_{k=1}^{\infty} \otimes \{d_k \}_{k=1}^{\infty} \in N(V).$ The series $\displaystyle \sum_{k=1}^{\infty} c_k f_{k+1} \otimes \sum_{k=1}^{\infty} d_k g_{k+1}$ converges by what is already proved, and furthermore

\begin{align*}
    \sum_{k=1}^{\infty} c_k f_{k+1} \otimes \sum_{k=1}^{\infty} d_k g_{k+1} &= \sum_{k=1}^{\infty} c_k T_1 f_k \otimes \sum_{k=1}^{\infty} d_k T_2 g_k \\
    &= \left( \widetilde{T_1} \otimes \widetilde{T_2} \right) \left( \sum_{k=1}^{\infty} c_k f_k \otimes \sum_{k=1}^{\infty} d_k g_k \right)\\
    &=0;
\end{align*}
this shows that $ \mathcal{T}_1 \{c_k \}_{k=1}^{\infty} \otimes \mathcal{T}_2 \{d_k \}_{k=1}^{\infty} \in \mathcal{N}_V.$\par

Finally, for every $k \in \mathbb{N},$
\begin{align*}
    \left| \left| f_{k+1} \otimes g_{k+1} \right| \right| &= \left| \left| T_1 f_k \otimes T_2 g_k \right| \right|\\
    &= \left| \left| (T_1 \otimes T_2 )(f_k \otimes g_k) \right| \right|\\
    & \leq \left| \left| T_1 \otimes T_2 \right| \right| \left| \left| f_k \otimes g_k \right| \right|,
\end{align*}

and thus $\left\{ \frac{ || f_{k+1} \otimes g_{k+1} ||}{|| f_k \otimes g_k ||} \right\}_{k=1}^{\infty} \in l^{\infty}.$\par
(2) Assume firstly that $T_1 \otimes T_2$ is bounded for every $\{ c_k \}_{k=1}^{\infty} \otimes \{ d_k \}_{k=1}^{\infty} \in \mathcal{D}(V)$. By (1) that $\mathcal{T}_1 \{c_k \}_{k=1}^{\infty} \otimes \mathcal{T}_2 \{d_k \}_{k=1}^{\infty} \in \mathcal{D}(V);$ furthermore,
\begin{align*}
    \left| \left| V ( \mathcal{T}_1 \{ c_k \}_{k=1}^{\infty} \otimes \mathcal{T}_2 \{d_k \}_{k=1}^{\infty} ) \right| \right| &= \left| \left| \sum_{k=1}^{\infty} c_k f_{k+1} \otimes \sum_{k=1}^{\infty} d_k g_{k+1} \right| \right|\\
    &= \left| \left| \sum_{k=1}^{\infty} c_k \widetilde{T_1}f_k \otimes \sum_{k=1}^{\infty} d_k \widetilde{T_2} g_k \right| \right|\\
    &= \left| \left| (\widetilde{T_1} \otimes \widetilde{T_2})\left( \sum_{k=1}^{\infty} c_k f_k \otimes \sum_{k=1}^{\infty} d_k g_k \right) \right| \right|\\
    & \leq \left| \left| \widetilde{T_1} \otimes \widetilde{T_2} \right| \right| \left| \left| V( \{c_k \}_{k=1}^{\infty} \otimes \{ d_k \}_{k=1}^{\infty} ) \right| \right|.
\end{align*}

Conversely, assume that there is a constant $K> 0$ so that $\left| \left| V( \mathcal{T}_1 \{c_k \}_{k=1}^{\infty} \otimes \mathcal{T}_2 \{d_k \}_{k=1}^{\infty} ) \right| \right| \leq K \left| \left| V ( \{c_k \}_{k=1}^{\infty} \otimes \{ d_k \}_{k=1}^{\infty} ) \right| \right|$ for all finite sequences.

Take an arbitrary $f \otimes g \in span \{ f_k \}_{k=1}^{\infty} \otimes span \{g_k \}_{k=1}^{\infty}$.i.e., $ \displaystyle  f \otimes g = \sum_{k=1}^{N} c_k f_k \otimes \sum_{k=1}^{N} d_k g_k$ for some $N \in \mathbb{N}$ and some $c_1, c_2,......., c_N, d_1, d_2, ......,d_N \in \mathbb{C};$ letting $c_k=0 , d_k =0$ for $k > N$ we have that 
\begin{align*}
    \left| \left| T_1 f \otimes T_2 g \right| \right| &= \left| \left| \sum_{k=1}^{\infty} c_k f_{k+1} \otimes \sum_{k=1}^{\infty} d_k g_{k+1} \right| \right|\\
     &= \left| \left| V ( \mathcal{T}_1 \{c_k \}_{k=1}^{\infty} \otimes \mathcal{T}_2 \{d_k \}_{k=1}^{\infty} ) \right| \right|\\
    & \leq K \left| \left| V( \{ c_k \}_{k=1}^{\infty} \otimes \{ d_k \}_{k=1}^{\infty}) \right| \right| \\
    &\leq K \left| \left| f \otimes g \right| \right|,
\end{align*}
i.e., $T_1 \otimes T_2$ is bounded.

\section{Frames and the Carleson condition in the tensor product of Hardy Spaces}
In this section, we consider a class of frames that can be represented via bounded operators $ T_1 \otimes T_2$. We extended the construction that first appeared in Corollary 3.17 in \cite{aldroubi2017dynamical} in terms of the tensor product. Our purpose is to extend the result given by Shapiro and Shields \cite{shapiro1961some} in the tensor product of Hardy spaces. First, we will discuss the Carleson condition on sequences $\{ \lambda_k \}_{k=1}^{\infty}$  and $ \{ \gamma_l \}_{l=1}^{\infty} $ in the open unit discs $\mathbb{D}_1$ and $\mathbb{D}_2$ respectively.
\subsection{The Carleson Condition  and tensor product}
Let $\mathbb{D}_1$ and $\mathbb{D}_2$ denote the open unit discs in the complex plane. The Hardy space $H_{1}^{2}(\mathbb{D}_1) \otimes H_{2}^{2}(\mathbb{D}_2)$ is defined by 
\newpage
 \[ H_{1}^{2}(\mathbb{D}_1) \otimes H_{2}^{2}(\mathbb{D}_2) = \biggl\{ f \otimes g : \mathbb{D}_1 \otimes \mathbb{D}_2 \to \mathbb{C} \big| f(z) \otimes g(\omega) = \sum_{n=0}^{\infty} a_n z^n \otimes \sum_{m=0}^{\infty} b_m \omega^m, \]
\[  \{a_n \}_{n=0}^{\infty} \otimes \{b_m \}_{m=0}^{\infty} \in l^2 ( \mathbb{N}_0 ) \otimes l^2 ( \mathbb{N}_0 ) \biggr\}.\]

\noindent The Hardy space $H_{1}^{2}(\mathbb{D}_1) \otimes H_{2}^{2}(\mathbb{D}_2) $ is a Hilbert space; given $w_1, w_2 \in  H_{1}^{2}(\mathbb{D}_1) \otimes H_{2}^{2}(\mathbb{D}_2)$ where $w_1 = f_1 \otimes g_1, w_2 = f_2 \otimes g_2 $ and $f_1 , f_2 \in H_{1}^{2}(\mathbb{D}_1)$ and $ g_1 , g_2 \in H_{2}^{2} (\mathbb{D}_2)$ 
\[ w_1 = \sum_{n=0}^{\infty} a_n z^n \otimes \sum_{m=0}^{\infty} b_m \omega^m , w_2 = \sum_{n=0}^{\infty} c_nz^n \otimes \sum_{m=0}^{\infty} d_m \omega^m \]
The inner product is defined by,

\begin{align*}
\displaystyle \langle w_1, w_2 \rangle &= \langle \sum_{n=0}^{\infty} a_n z^n \otimes \sum_{m=0}^{\infty} b_m \omega^m , \sum_{n=0}^{\infty} c_nz^n \otimes \sum_{m=0}^{\infty} d_m \omega^m \rangle \\
 &= \langle \sum_{n=0}^{\infty} a_n z^n, \sum_{n=0}^{\infty} c_n z^n \rangle \langle \sum_{m=0}^{\infty}b_m \omega^m , \sum_{m=0}^{\infty} d_m \omega^m \rangle\\
 &= \sum_{n=0}^{\infty} a_n \bar{c_n} \cdot \sum_{m=0}^{\infty} b_m \bar{d_m}
\end{align*} 

Note that $\{ z^n\}_{n=0}^{\infty} \otimes \{ \omega^m \}_{m=0}^{\infty}$ is an orthonormal basis for $H_{1}^{2}(\mathbb{D}_1) \otimes H_{2}^{2}(\mathbb{D}_2)$; denoting the canonical basis for $l^2(\mathbb{N}) \otimes l^2(\mathbb{N})$ by $\{ \delta_n \}_{n=1}^{\infty} \otimes \{ \sigma_m \}_{m=1}^{\infty},$ the operator $\theta : H_{1}^{2}(\mathbb{D}_1) \otimes H_{2}^{2}(\mathbb{D}_2) \to l^2(\mathbb{N}) \otimes l^2(\mathbb{N})$ defined by $\theta (z^n \otimes \omega^m) = \delta_{n+1} \otimes \sigma_{n+1}$ for $n,m = 0,1,2,3,.......$ is a unitary operator from   $H_{1}^{2}(\mathbb{D}_1) \otimes H_{2}^{2}(\mathbb{D}_2)$ onto $l^2(\mathbb{N}) \otimes l^2(\mathbb{N})$ and its preserve the norm and inner product of vectors.

\begin{definition} \rm
A sequence $\{\lambda_k \}_{k=1}^{\infty} \subset \mathbb{D}_1, \{ \gamma_l \}_{l=1}^{\infty} \subset \mathbb{D}_2$ then $\{\lambda_k \}_{k=1}^{\infty} \otimes \{ \gamma_l \}_{l=1}^{\infty}$ satisfies the Carleson condition if 
\begin{equation}  
 \displaystyle \inf_{n,m \in \mathbb{N}} \prod_{(k,l) \neq (n,m) } \frac{\mid \lambda_k \gamma_l - \lambda_n \gamma_m \mid}{\mid 1 - \overline{\lambda_k \gamma_l} \lambda_n \gamma_m \mid } > 0.\end{equation}

\end{definition}

\noindent    For a given sequence $\Lambda = \{ \lambda_{k} \}_{k=1}^{\infty} \otimes \{ \gamma_{l}\}_{l=1}^{\infty} \subset \mathbb{D}_1 \otimes \mathbb{D}_2, $ define the sequence-valued operator $\varphi_{\Lambda}$ by

    \footnotesize\begin{equation}
        \varphi_{\Lambda}(f \otimes g) = \{ f(\lambda_k) \sqrt{ 1 - | \lambda_{k} |^2 }\}_{k=1}^{\infty} \otimes \{ g(\gamma_l )\sqrt{ 1-|\gamma_l |} \}_{l=1}^{\infty} , f\otimes g \in H_{1}^{2}(\mathbb{D}_1) \otimes H_{2}^{2}(\mathbb{D}_2).\end{equation}

        \noindent Since the sequence in equation (2) does not necessarily belong to $l^2(\mathbb{N}) \otimes l^2( \mathbb{N}).$ The following result is motivated by the result given by Shapiro and Shields in \cite{shapiro1961some} and the Theorem 9.1 of \cite{duren1970theory}.

\begin{theorem}\rm
    A sequence $\{ \lambda_k \}_{k=1}^{\infty} \otimes \{ \gamma_l \}_{l=1}^{\infty} \subset \mathbb{D}_1 \otimes \mathbb{D}_2$ satisfies the Carleson condition.i.e.,\\
    $\displaystyle \prod_{k,l = 1, (l \neq m) }^{\infty} \frac{\mid \lambda_k \gamma_l - \lambda_n \gamma_m \mid}{\mid 1 - \overline{\lambda_k \gamma_l} \lambda_n \gamma_m \mid } > 0$ if and only if the interpolation problem $f(\lambda_k) g(\gamma_l) = z_k \omega_l, f\otimes g $ are bounded analytic in $\mathbb{D}_1 \otimes \mathbb{D}_2$ is solvable for arbitrary $z_k \omega_l$. i.e., $\varphi_{\Lambda}( H_{1}^{2}(\mathbb{D}_1) \otimes H_{2}^{2}(\mathbb{D}_2)) = l^2(\mathbb{N}) \otimes l^2(\mathbb{N})$ and in affirmative case, $\varphi_{\Lambda}$ is bounded.
\end{theorem}

\begin{theorem}\rm
Let $\{\lambda_k \}_{k=1}^{\infty} \subset \mathbb{D}_1$ and $\{ \gamma_l \}_{l=1}^{\infty} \subset \mathbb{D}_2$ be a sequence of distinct numbers and $\{ \lambda_k \}_{k=1}^{\infty} \cap \{ \gamma_l \}_{l=1}^{\infty} = \phi$ then $\{ \lambda_k \}_{k=1}^{\infty} \otimes \{ \gamma_l \}_{l=1}^{\infty} \subset \mathbb{D}_1 \otimes \mathbb{D}_2$ is also sequence of distinct numbers. If $\exists \; c \in (0,1)$ such that 
\begin{equation}    
\dfrac{1 - \mid \lambda_k \gamma_{l+1}\mid}{1 - \mid \lambda_k \gamma_l \mid} \leq c <1, \forall k,l \in \mathbb{N}  \end{equation}

then $\{ \lambda_k \}_{k=1}^{\infty} \otimes \{ \lambda_l \}_{l=1}^{\infty}$ satisfies the carleson condition. If $ \{ \lambda_k \}_{k=1}^{\infty} \otimes \{ \gamma_l \}_{l=1}^{\infty} $ is positive and increasing then condition (3) is also necessary for $ \{\lambda_k \}_{k=1}^{\infty} \otimes \{ \gamma_l \}_{l=1}^{\infty}$ to satisfy the Carleson condition.
\end{theorem}
\textbf{Proof :} Condition (3) implies that 
\[ 1 - \mid \lambda_k \gamma_{l+1} \mid \leq c ( 1- \lambda_k \gamma_l )\]
\begin{equation}
    1 - \mid \lambda_k \gamma_l \mid \leq c^{l-m}( 1 - \lambda_n \gamma_m), l >m, k \geq n \; \text{and}\;  k,m = 1,2,3.... 
\end{equation}

    In particular 
    $\displaystyle \sum_{k} \sum_{l} ( 1- \mid \lambda_k \gamma_l \mid) < \infty$

    it follows from (4) that for $ l > m$

    \[ \mid \lambda_k \gamma_l \mid - \mid \lambda_n \gamma_m \mid \geq ( 1- c^{l-m})( 1- \mid \lambda_n \gamma_m \mid)\]

    and \begin{align*}
        1 - \mid \lambda_k \gamma_l \lambda_n \gamma_m \mid &= 1 - \mid \lambda_k \gamma_l \mid + \mid \lambda_k \gamma_l \mid ( 1- \mid \lambda_n \gamma_m \mid )\\
        & \leq ( 1+ c^{l-m})(1 - \mid \lambda_n \gamma_m \mid )
    \end{align*}

    Hence by the lemma 
    \begin{align*}
        \left| \dfrac{ \lambda_n \gamma_m - \lambda_k \gamma_l }{1 - \overline{ \lambda_k \gamma_l } \lambda_n \gamma_m}\right| & \geq \dfrac{ \mid \lambda_k \gamma_l \mid - \mid \lambda_n \gamma_m \mid}{1 - \mid \lambda_k \gamma_l \lambda_n \gamma_m \mid }\\
        &\geq \dfrac{ 1- c^{l-m}}{1 + c^{l-m}},\; \; \; \; \; \;  l > m
    \end{align*}

    for $ l < m$ this inequality takes the form 

\[ \left| \dfrac{ \lambda_m \gamma_m - \lambda_k \gamma_l}{ 1- \overline{\lambda_k \gamma_l}\lambda_n \gamma_m } \right| \geq \dfrac{ 1 - c^{m-l}}{1 + c^{m-l}}\]

Consequently,

\[ \prod_{k,l=1, (l \neq m)}^{\infty} \left| \dfrac{\lambda_n \gamma_m - \lambda_k \gamma_l}{1 - \overline{\lambda_k \gamma_l}\lambda_n \gamma_m}\right| \geq \prod_{N=1}^{\infty} \left( \dfrac{1 - c^{N}}{1 + c^{N}} \right)^2 > 0\]

which shows that $ \{ \lambda_k \}_{k=1}^{\infty} \otimes \{ \gamma_l \}_{l=1}^{\infty}$ satisfies the carleson condition. Now suppose $ 0 \leq \lambda_1 \gamma_1 < \lambda_1 \gamma_2 < .........< \lambda_2 \gamma_1 < \lambda_2 \gamma_2 < ......$ and 

\[ \prod_{k,l=1, ( l \neq m)}^{\infty} \left| \dfrac{\lambda_n \gamma_m - \lambda_k \gamma_l}{1 - \lambda_k \gamma_l \lambda_n \gamma_m } \right| \geq \delta > 0\]

Then, \[ \lambda_k \gamma_{l+1} - \lambda_k \gamma_l \geq \delta ( 1 - \lambda_k \gamma_l \lambda_k \gamma_{l+1}),\; \; \; \; \; \; \; \; \; \; \; \; \; \; \; k ,l = 1,2,.....\]

So that,

\[ 1 - \lambda_k \gamma_{l+1} \leq 1 - \dfrac{ \delta + \lambda_k \gamma_l}{1 + \delta \lambda_k \gamma_l } \leq ( 1- \delta)(1 - \lambda_k \gamma_l)\]

Thus $ \{ \lambda_k \}_{k=1}^{\infty} \otimes \{ \gamma_l \}_{l=1}^{\infty}$ satisfies (3) as claimed.


    \par

    \subsection{Frame Properties and the Carleson Condition in terms of tensor product}

    In this section, we have extended the result given Theorem 3.7 in \cite{christensen2018operator}, which yields the construction of a class of operators $T_1 \otimes T_2: l^2(\mathbb{N}) \otimes l^2(\mathbb{N}) \to l^2(\mathbb{N}) \otimes l^2(\mathbb{N})  $ for which $\{T_{1}^{n}h_1 \}_{n=0}^{\infty} \otimes \{T_{2}^{m} h_2 \}_{m=0}^{\infty}$ is a frame for $l^2(\mathbb(N)) \otimes l^2(\mathbb(N))$ for certain sequences $h_1 \otimes h_2 \in l^2(\mathbb(N)) \otimes l^2(\mathbb(N)). $

    Consider a sequence $\{\lambda_k \}_{k=1}^{\infty} \subset \mathbb{D}_1$ and $\{ \gamma_{l} \}_{l=1}^{\infty} \subset \mathbb{D}_2$, where $\mathbb{D}_1, \mathbb{D}_2$ are unit discs in arbitrary complex planes and assume that $ \displaystyle \{ \sqrt{1 - \mid \lambda_{k} \mid^2} \}_{k=1}^{\infty} \in l^2( \mathbb{N}) $ and $\{ \sqrt{1 - \mid \gamma_l \mid^2} \}_{l=1}^{\infty} \in l^2( \mathbb{N})$. Given any Hilbert spaces $\mathcal{H}_1$ and $\mathcal{H}_2$, choose an orthonormal basis $\{ e_k \}_{k=1}^{\infty} $ and $\{ d_l \}_{l=1}^{\infty}$ and consider the bounded linear operators $T_1 : \mathcal{H}_1 \to \mathcal{H}_1, T_2: \mathcal{H}_2 \to \mathcal{H}_2 $ for which $T_1e_k = \lambda_k e_k, T_2d_l = \gamma_l d_l.$ Let $\displaystyle h_1 = \sum_{k=1}^{\infty} \sqrt{ 1 - \mid \lambda_k \mid^2}e_k, h_2 = \sum_{l=1}^{\infty} \sqrt{ 1 - \mid \gamma_l \mid^2}d_l $ and consider the iterated systems
\begin{equation}  
\{ T_{1}^{n} h_1 \}_{n=0}^{\infty} = \{ \sum_{k=1}^{\infty} \lambda_{k}^{n} \sqrt{ 1 - \mid \lambda_k \mid^2} e_k \}_{n=0}^{\infty}
\end{equation}

\begin{equation}
    \{ T_{2}^{m} h_2 \}_{m=0}^{\infty} = \{ \sum_{l=1}^{\infty} \gamma_{l}^{m} \sqrt{ 1 - \mid \gamma_l \mid^2} d_l \}_{m=0}^{\infty} 
    \end{equation}
We will now state the following result.

\begin{theorem}\rm
Let $\{ \lambda_k \}_{k=1}^{\infty} \subset \mathbb{D}_1, \{ \gamma_l \}_{l=1}^{\infty} \subset \mathbb{D}_2 $ and assume that $\{ \sqrt{1 - \mid \lambda_k \mid^2} \}_{k=1}^{\infty} \in l^2( \mathbb{N}), \\ \{ \sqrt{1 - \mid \gamma_l \mid^2} \}_{l=1}^{\infty} \in l^2( \mathbb{N})$ Then the sequence $\{ T_{1}^{n} h_1 \}_{n=0}^{\infty} \otimes \{ T_{2}^{m} h_2 \}_{m=0}^{\infty}$ defined by $(5)$ and $(6)$ is a frame for $\mathcal{H}_1 \otimes \mathcal{H}_2$ if and only if $\{ \lambda_k \}_{k=1}^{\infty} \otimes \{ \gamma_l \}_{l=1}^{\infty}$ satisfies the Carleson condition.

\end{theorem}

\textbf{Proof :}
Define the synthesis operator $ V : l^2(\mathbb{N}_0) \otimes l^2(\mathbb{N}_0) \to \mathcal{H}_1 \otimes \mathcal{H}_2$ by $V\{  \{a_n \}_{n=0}^{\infty} \otimes \{b_m \}_{0}^{\infty} \} =\displaystyle \sum_{n=0}^{\infty} a_n T_{1}^n h_1 \otimes  \sum_{m=0}^{\infty} b_m T_{2}^{m} h_2$. The sequence $\{T_{1}^{n} h_1 \}_{n=0}^{\infty} \otimes \{T_{2}^{m} h_2 \}_{m=0}^{\infty} $ is a frame for $\{ \mathcal{ H}_1\} \otimes \{\mathcal{H}_2 \}$ if and only if the operator $V$ is well defined and surjective.

       First assume that $ \{T_{1}^{n} h_1 \}_{n=0}^{\infty} \otimes \{T_{2}^{m} h_2 \}_{m=0}^{\infty}$ is a frame for $\mathcal{H}_1 \otimes \mathcal{H}_2$. Take an arbitrary sequence $ \{c_j \}_{j=1}^{\infty} \otimes \{ f_r \}_{r=1}^{\infty} \in l^2( \mathbb{N} ) \otimes l^2( \mathbb{N})$ The surjectivity of $V$ implies that there exists $\{a_n\}_{n=0}^{\infty} \otimes \{b_m \}_{m=0}^{\infty} \in l^2( \mathbb{N}_0) \otimes l^2(\mathbb{N}_0)$ such that $ \displaystyle \sum_{n=0}^{\infty} a_n T_{1}^{n} h_1 \otimes \sum_{m=0}^{\infty} b_mT_{2}^{m} h_2 = \sum_{j=1}^{\infty} c_j e_j \otimes \sum_{r=1}^{\infty} f_r d_r.$ It follows that for each $k,l \in \mathbb{N},$

\begin{equation}    
\begin{split}
     c_k f_l &= \langle \sum_{j=1}^{\infty} c_j e_j \otimes \sum_{r=1}^{\infty} f_r d_r, e_k \otimes d_l \rangle\\
&= \langle \sum_{j=1}^{\infty}c_j e_j, e_k \rangle \langle \sum_{r=1}^{\infty} f_rd_r, d_l \rangle \\
&= \sum_{n=0}^{\infty} a_n \langle T_{1}^{n}h_1, e_k \rangle \sum_{m=0}^{\infty} b_m \langle T_{2}^{m}h_2, d_l \rangle \\
&= \sum_{n=0}^{\infty} a_n \lambda_{k}^{n}\sqrt{1- | \lambda_k|^2} \cdot \sum_{m=0}^{\infty} b_m \gamma_{l}^{m} \sqrt{1 - | \gamma_l |^2},
\end{split}
\end{equation}
where $c_k = \displaystyle \langle \sum_{j=1}^{\infty} c_j e_j, e_k \rangle = \sum_{n=0}^{\infty} a_n \left< T_{1}^{n}h_1, e_k \right> = \sum_{n=0}^{\infty} a_n \lambda_{k}^{n} \sqrt{ 1- \mid \lambda_k \mid^2 }$ and similarly \\$ \displaystyle f_l = \sum_{m=0}^{\infty} b_m \gamma_{l}^{m} \sqrt{ 1 - \mid \gamma_l \mid^2}.$\\
\par
 
 Defining $f\in H_{1}^{2}(\mathbb{D}_1)$ and $g \in H_{2}^{2}( \mathbb{D}_2)$ by 

 \[ f(z) = \sum_{n=0}^{\infty} a_n z^n,  g( \omega) = \sum_{m=0}^{\infty} b_m \omega^m. \] Equation $(7)$ turns into 
 \[ f( \lambda_k)\sqrt{1- |\lambda_k |^2} \cdot g(\gamma_l) \sqrt{1 - | \gamma_l |^2} = c_kf_l. \] Formulated in terms of the operator $\varphi_{\Lambda}$ in $(2)$, this means that $l^2(\mathbb{N} ) \otimes l^2( \mathbb{N}) \subseteq \varphi_{\Lambda}( H_{1}^{2}(\mathbb{D}_1) \otimes H_{2}^{2}(\mathbb{D}_2))$.\\ \par
 
 On the other hand, take an arbitrary $f \otimes g \in H_{1}^{2}(\mathbb{D}_1) \otimes H_{2}^{2}(\mathbb{D}_2)$ and choose $\{ a_n \}_{n=0}^{\infty} \otimes \{b_m \}_{m=0}^{\infty} \in l^2( \mathbb{N}_0) \otimes l^2( \mathbb{N}_0)$ such that, $\displaystyle f(z) = \sum_{n=0}^{\infty} a_n z^n, g(\omega) = \sum_{m=0}^{\infty} b_m \omega^m$ for every $k, l \in \mathbb{N}$ we have 

\footnotesize\begin{align*}
    \langle V\{ \{a_n\}_{n=0}^{\infty} \otimes \{b_m \}_{m=0}^{\infty} \} , e_k \otimes d_l \rangle &= \langle \sum_{n=0}^{\infty} a_n T_{1}^{n} h_1 \otimes \sum_{m=0}^{\infty} b_m T_{2}^{m} h_2, e_k \otimes d_l \rangle\\
    &= \langle \sum_{n=0}^{\infty} a_n \sum_{j=1}^{\infty} \lambda_{j}^{n} \sqrt{ 1 - |\lambda_j|^2}e_j \otimes \sum_{m=0}^{\infty} b_m \sum_{r=1}^{\infty} \gamma_{r}^{m} \sqrt{1- | \gamma_r |^{2}} d_r\\
    & \hspace{5cm}, e_k \otimes d_l \rangle\\
&= \sum_{n=0}^{\infty} a_n \lambda_{k}^{n} \sqrt{ 1- |\lambda_k |^2} \cdot \sum_{m=0}^{\infty} b_m \gamma_{l}^{m}\sqrt{ 1- | \gamma_{l} |^2} 
\end{align*}

Therefore, \[ \varphi_{\Lambda} (f \otimes g) =  \{ \langle V\{ \{a_n\}_{n=0}^{\infty} \otimes \{b_m \}_{m=0}^{\infty} \} , e_k \otimes d_l \rangle \}_{k,l =1}^{\infty} \in l^2( \mathbb{N}) \otimes l^2( \mathbb{N}) \]

\[\implies \varphi_{\Lambda}( H_{1}^{2}(\mathbb{D}_1) \otimes H_{2}^{2}(\mathbb{D}_2)) \subseteq l^2(\mathbb{N}) \otimes l^2(\mathbb{N})\]

Hence, \[ \varphi_{\Lambda}( H_{1}^{2}(\mathbb{D}_1) \otimes H_{2}^{2}(\mathbb{D}_2)) = l^2(\mathbb{N}) \otimes l^2(\mathbb{N}),\] which by [Theorem 3] implies that $\{\lambda_k \}_{k=1}^{\infty} \otimes \{\gamma_l\}_{l=1}^{\infty}$ satisfies the Carleson condition.\\

        Conversely, assume that the sequence $\{ \lambda_j \}_{j=1}^{\infty} \otimes \{ \gamma_l \}_{l=1}^{\infty} \subset \mathbb{D}_1 \otimes \mathbb{D}_2$ satisfies the Carleson condition. We first show that $\displaystyle \sum_{n=0}^{\infty} a_n T_{1}^{n} h_1 \otimes \sum_{m=0}^{\infty} b_m T_{2}^{m} h_2$ is convergent for all $\{a_n \}_{n=0}^{\infty} \otimes \{ b_m \}_{m=0}^{\infty} \in l^2( \mathbb{N}_0) \otimes l^2( \mathbb{N}_0)$. Consider the corresponding $f \in H_{1}^{2}(\mathbb{D}_1)$ and $g \in H_{2}^{2}(\mathbb{D}_2)$ determined by $\displaystyle f(z) = \sum_{n=0}^{\infty}a_nz^n, g(\omega) = \sum_{m=0}^{\infty}b_m\omega^m$. By [Theorem 3] we know that $\varphi_{\Lambda}( H_{1}^{2}(\mathbb{D}_1) \otimes H_{2}^{2}(\mathbb{D}_2)) = l^2( \mathbb{N}) \otimes l^2( \mathbb{N})$. So $\displaystyle \{ f(\lambda_k) \sqrt{1 - |\lambda_k |^2} \}_{k=1}^{\infty} \otimes \{ g( \gamma_l) \sqrt{ 1 - |\gamma_l |^2} \}_{l=1}^{\infty} \in l^2(\mathbb{N}) \otimes l^2( \mathbb{N})$. Now for $N,M \in \mathbb{N}$, consider the truncated sequences $\displaystyle \{ a_n \}_{n=0}^{N}, \{ b_m \}_{m=0}^{M} $ and the associated functions $\displaystyle f_N \in H_{1}^{2}(\mathbb{D}_1), g_{M} \in H_{2}^{2}( \mathbb{D}_2)$ ,respectively. Given by $\displaystyle f_N(z) = \sum_{n=0}^{N}a_nz^n, g_M(\omega) = \sum_{m=0}^{M} b_m \omega^m.$\\
        Again $\{ f(\lambda_k) \sqrt{1 - |\lambda_k |^2} \}_{k=1}^{\infty} \otimes \{ g( \gamma_l) \sqrt{ 1 - |\gamma_l |^2} \}_{l=1}^{\infty} \in l^2(\mathbb{N}) \otimes l^2( \mathbb{N})$, and since $\varphi_\Lambda: H_{1}^{2}(\mathbb{D}_1) \otimes H_{2}^{2}(\mathbb{D}_2) \to l^2( \mathbb{N}) \otimes l^2( \mathbb{N})$ is bounded by [Theorem 3], there is a constant $C>0$ such that 

\footnotesize\begin{align*}     
         || \varphi_{\Lambda}( f \otimes g ) - \varphi_{\Lambda}(f_N \otimes g_M) ||^2 &\leq C || f \otimes g - f_N \otimes g_M ||^2 \\
&= C|| \sum_{n=0}^{\infty} a_n z^n \sum_{m=0}^{\infty} b_m \omega^m - \sum_{n=0}^{N} a_n z^n \sum_{m=0}^{M} b_m \omega^m ||^2\\
\footnotesize&= C||\left( \sum_{n=0}^{N} a_n z^n + \sum_{n=N+1}^{\infty} a_n z^n \right) \left( \sum_{m=0}^{M} b_m \omega^m + \sum_{m=M+1}^{\infty} b_m \omega^m \right)\\
&\; \; \; \; \; \; \hspace{2cm}  -\sum_{n=0}^{N} a_n z^n \sum_{m=0}^{M} b_m \omega^m ||^2\\
\footnotesize&= C || \sum_{n=0}^{N}a_n z^n \sum_{m=0}^{M} b_m \omega^m + \sum_{n=0}^{N} a_n z^n \sum_{m=M+1}^{\infty} b_m \omega^m \\
 & \hspace{2cm} + \sum_{n=N+1}^{\infty} a_n z^n \sum_{m=0}^{M} b_m \omega^m + \sum_{n=N+1}^{\infty} a_n z^n \sum_{m= M+1}^{\infty} b_m \omega^m\\
& \hspace{3cm}- \sum_{n=0}^{N} a_n z^n \sum_{m=0}^{M} b_m \omega^m ||^2 \\
&= C ||\sum_{n=0}^{N} a_n z^n \sum_{m=M+1}^{\infty} b_m \omega^m  + \sum_{n=N+1}^{\infty} a_n z^n \sum_{m=0}^{M} b_m \omega^m\\
\displaybreak& \hspace{2cm} + \sum_{n=N+1}^{\infty} a_n z^n \sum_{m= M+1}^{\infty} b_m \omega^m ||^2 \\
&= C \sum_{n=0}^{N} | a_n |^2 \sum_{m = M+1}^{\infty} |b_m |^2 + C\sum_{n=N+1}^{\infty} |a_n|^2 \sum_{m=0}^{M} | b_m |^2\\
& \hspace{3cm}+ C\sum_{n=N+1}^{\infty} | a_n |^2 \sum_{m= M+1}^{\infty} |b_m |^2  \to 0, N,M \to \infty\\
 \end{align*}

It follows that 
\begin{align*}
 \sum_{n=0}^{N} a_n T_{1}^{n}h_1 \otimes \sum_{m=0}^{M} b_mT_{2}^{m}h_2 &= \sum_{n=0}^{N}a_n \sum_{k=1}^{\infty} \sqrt{1 - \mid \lambda_k \mid^2} \lambda_{k}^{n} e_k \otimes \sum_{m=0}^{M} b_m \sum_{l=1}^{\infty} \sqrt{1 - \mid \gamma_l \mid^2}\gamma_{l}^{n} d_l\\
&= \sum_{k=1}^{\infty} \sqrt{1 - \mid \lambda_k \mid^k}\sum_{n=0}^{N} a_n \lambda_{k}^{n} e_k \cdot \sum_{l=1}^{\infty} \sqrt{1 - \mid \gamma_l \mid^2}\sum_{m=0}^{M}b_m \gamma_{l}^{m} d_l\\
 &= \sum_{k=1}^{\infty} \sqrt{1 - \mid \lambda_{k} \mid^2} f_{N}(\lambda_k) e_k \cdot \sum_{l=1}^{\infty} \sqrt{ 1 - \mid \lambda_{l} \mid^2} g_{M}( \gamma_l ) d_l\\
& \to \sum_{k=1}^{\infty} f( \lambda_k ) \sqrt{1 - \mid \lambda_k \mid^2} e_k \sum_{l=1}^{\infty} g(\gamma_l) \sqrt{1 - \mid \gamma_l \mid^2} d_l \; \; \; \; \text{as } N, M \to \infty\\    
\end{align*}

This proves that $\displaystyle \sum_{n=0}^{\infty} a_n T_{1}^n h_1 \otimes \sum_{m=0}^{\infty} b_m T_{2}^{m} h_2$ is convergent, and thus $V$ is well defined from $l^2(\mathbb{N}_0) \otimes l^2( \mathbb{N}_0)$ into $\mathcal{H}_1 \otimes \mathcal{H}_2$. In order to prove that $\{ T_{1}^n h_1 \}_{n=0}^{\infty} \otimes \{ T_{2}^{m}h_2 \}_{m=0}^{\infty}$ is a frame, it is enough to show that the synthesis operator $V: l^2(\mathbb{N}_0) \otimes l^2(\mathbb{N}_0) \to \mathcal{H}_1 \otimes \mathcal{H}_2$ is surjective. Let $x \otimes y \in \mathcal{H}_1 \otimes \mathcal{H}_2$.By [Theorem 3] there is an $ f \otimes g \in H_{1}^{2}(\mathbb{D}_1) \otimes H_{2}^{2}(\mathbb{D}_2)$ such that $f(\lambda_k )\sqrt{1 - \mid \lambda_k \mid^2} \otimes g( \gamma_l)\sqrt{ 1 -\mid \gamma_l \mid^2} = \langle x \otimes y, e_k \otimes d_l \rangle $ for all $k,l \in \mathbb{N}$. Choose $\{ a_n \}_{n=0}^{\infty} \otimes \{b_m \}_{m=0}^{\infty} \in l^2( \mathbb{N}_0) \otimes l^2( \mathbb{N}_0)$, such that $\displaystyle f(z) = \sum_{n=0}^{\infty} a_n z^n, g(\omega ) = \sum_{m=0}^{\infty} b_m \omega^m$. Then for each $k,l \in \mathbb{N},$

\begin{align*}
    \langle V (\{ a_n \}_{n=0}^{\infty} \otimes \{b_m \}_{m=0}^{\infty}), e_k \otimes d_l \rangle &= \langle \sum_{n=0}^{\infty} a_n T_{1}^n h_1 \otimes \sum_{m=0}^{\infty} b_m T_{2}^m h_2, e_k \otimes d_l \rangle \\
    &= \langle \sum_{n=0}^{\infty} a_n T_{1}^n h_1, e_k \rangle \langle \sum_{m=0}^{\infty} b_m T_{2}^{m} h_2, d_l \rangle \\
    &= \sum_{n=0}^{\infty} a_n \langle \sum_{j=1}^{\infty} \lambda_{j}^{n} \sqrt{1 - \mid \lambda_j \mid^2 } e_j , e_k \rangle \\
    \displaybreak & \hspace{2cm} \cdot\sum_{m=0}^{\infty} b_m \langle \sum_{r=1}^{\infty} \gamma_{r}^{m} \sqrt{ 1 - \mid \gamma_{r} \mid^2} d_r, d_l \rangle \\
    &= \sum_{n=0}^{\infty}a_n \lambda_{k}^{n} \sqrt{1 - \mid \lambda_k \mid^2} \sum_{m=0}^{\infty} b_m \gamma_{l}^{m} \sqrt{1 - \mid \gamma_l \mid^2 }\\
    &= f(\lambda_k) \sqrt{1 - \mid \lambda_k \mid^2} \cdot g( \gamma_l) \sqrt{ 1- \mid \gamma_l \mid^2}\\
     &= \langle x, e_k \rangle \langle y, d_l \rangle .\\
\end{align*}
 Therefore $V( \{a_n \}_{n=0}^{\infty} \otimes \{ b_m \}_{m=0}^{\infty} )= x\otimes y$ and thus $V$ is surjective , as desired.\\~\\

 \noindent {\bf Acknowledgement: } This work was supported by the Science and Engineering Research Board, DST, Govt. of India
 [CRG/2020/003170].\\
{\bf Declaration: }
 I confirm that all authors listed on the title page have contributed equally to the work, have read the manuscript,
  and agree to its submission.\\
 {\bf Conflict of Interest: } The authors declare that there is no
 conflict of interest.\\
 {\bf Data Availability:} Data availability or data sharing is not applicable to this article as no data sets were generated or analyzed during the
 current investigation.

\bibliographystyle{plain}
\bibliography{bibv}

\begin{thebibliography}{10}

\bibitem{aldroubi2017dynamical}
Akram Aldroubi, Carlos Cabrelli, Ursula Molter, and Sui Tang.
\newblock Dynamical sampling.
\newblock {\em Applied and Computational Harmonic Analysis}, 42(3):378--401,
  2017.

\bibitem{aldroubi2013dynamical}
Akram Aldroubi, Jacqueline Davis, and Ilya Krishtal.
\newblock Dynamical sampling: Time--space trade-off.
\newblock {\em Applied and Computational Harmonic Analysis}, 34(3):495--503,
  2013.

\bibitem{aldroubi2017dynamical1}
Akram Aldroubi and Armenak Petrosyan.
\newblock Dynamical sampling and systems from iterative actions of operators.
\newblock In {\em Frames and Other Bases in Abstract and Function Spaces: Novel
  Methods in Harmonic Analysis, Volume 1}, pages 15--26. Springer, 2017.

\bibitem{ashbrock2023dynamical}
Jonathan Ashbrock and Alexander~M Powell.
\newblock Dynamical dual frames with an application to quantization.
\newblock {\em Linear Algebra and its Applications}, 658:151--185, 2023.

\bibitem{bourouihiya2008tensor}
Abdelkrim Bourouihiya.
\newblock The tensor product of frames.
\newblock {\em Sampling theory in signal and Image processing}, 7:65--76, 2008.

\bibitem{cabrelli2020dynamical}
Carlos Cabrelli, Ursula Molter, Victoria Paternostro, and Friedrich Philipp.
\newblock Dynamical sampling on finite index sets.
\newblock {\em Journal d'analyse math{\'e}matique}, 140(2):637--667, 2020.

\bibitem{carando2009duality}
Daniel Carando and Silvia Lassalle.
\newblock Duality, reflexivity and atomic decompositions in banach spaces.
\newblock {\em Studia Math}, 191(1):67--80, 2009.

\bibitem{chen2017frame}
Qiuhui Chen, Pei Dang, and Tao Qian.
\newblock A frame theory of hardy spaces with the quaternionic and the clifford
  algebra settings.
\newblock {\em Advances in Applied Clifford Algebras}, 27:1073--1101, 2017.

\bibitem{christensen2017operator}
Ole Christensen and Marzieh Hasannasab.
\newblock Operator representations of frames: boundedness, duality, and
  stability.
\newblock {\em Integral Equations and Operator Theory}, 88:483--499, 2017.

\bibitem{christensen2019frame}
Ole Christensen and Marzieh Hasannasab.
\newblock Frame properties of systems arising via iterated actions of
  operators.
\newblock {\em Applied and Computational Harmonic Analysis}, 46(3):664--673,
  2019.

\bibitem{christensen2018dynamical}
Ole Christensen, Marzieh Hasannasab, and Ehsan Rashidi.
\newblock Dynamical sampling and frame representations with bounded operators.
\newblock {\em Journal of Mathematical Analysis and Applications},
  463(2):634--644, 2018.

\bibitem{christensen2018operator}
Ole Christensen, Marzieh Hasannasab, and Diana~T Stoeva.
\newblock Operator representations of sequences and dynamical sampling.
\newblock {\em Sampling Theory in Signal and Image Processing}, 17:29--42,
  2018.

\bibitem{duren1970theory}
Peter~L Duren.
\newblock {\em Theory of H p Spaces}.
\newblock Academic press, 1970.

\bibitem{khosravi2012frames}
Amir Khosravi and Mohammad~Sadegh Asgari.
\newblock Frames and bases in tensor product of hilbert spaces.
\newblock {\em arXiv preprint arXiv:1204.0096}, 2012.

\bibitem{martin2021continuous}
Roc{\'\i}o~D{\'\i}az Mart{\'\i}n, Ivan Medri, and Ursula Molter.
\newblock Continuous and discrete dynamical sampling.
\newblock {\em Journal of Mathematical Analysis and Applications},
  499(2):125060, 2021.

\bibitem{nemec2017quantitative}
Mike Nemec and Daniel Hoffmann.
\newblock Quantitative assessment of molecular dynamics sampling for flexible
  systems.
\newblock {\em Journal of chemical theory and computation}, 13(2):400--414,
  2017.

\bibitem{shapiro1961some}
Harold~S Shapiro and Allen~L Shields.
\newblock On some interpolation problems for analytic functions.
\newblock {\em American Journal of Mathematics}, 83(3):513--532, 1961.

\bibitem{speransky2019existence}
KS~Speransky and PA~Terekhin.
\newblock Existence of frames based on the szeg{\"o} kernel in the hardy space.
\newblock {\em Russian Mathematics}, 63:51--61, 2019.

\bibitem{stoeva2011characterization}
Diana~T Stoeva.
\newblock Characterization of atomic decompositions, banach frames, xd-frames,
  duals and synthesis-pseudo-duals, with application to hilbert frame theory.
\newblock {\em arXiv preprint arXiv:1108.6282}, 2011.

\bibitem{volkmer1995frames}
Hans Volkmer.
\newblock Frames of wavelets in hardy space.
\newblock {\em Analysis}, 15(4):405--422, 1995.

\end{thebibliography}

        \end{document}